# Frequency-Reactive Power Optimization Strategy of Grid-forming Offshore Wind Farm Using DRU-HVDC Transmission

Zhekai Li
*Key Laboratory of Control of Power Transmission and Conversion*
*Shanghai Jiao Tong University*
Shanghai, China
zhekai_li.work@sjtu.edu.cn

Kun Han
*Xuji Electric Co., Ltd.*
Xuchang China
hk20031009@163.com

Xu Cai
*Key Laboratory of Control of Power Transmission and Conversion*
*Shanghai Jiao Tong University*
Shanghai, China
xucai@sjtu.edu.cn

Renxin Yang
*Key Laboratory of Control of Power Transmission and Conversion*
*Shanghai Jiao Tong University*
Shanghai, China
frank_yang@sjtu.edu.cn

Haotian Yu
*Key Laboratory of Control of Power Transmission and Conversion*
*Shanghai Jiao Tong University*
Shanghai, China
yht56114@sjtu.edu.cn

Kepeng Xia
*Xuji Electric Co., Ltd.*
Xuchang China
283691206@qq.com

Lulu Liu
*Xuji Electric Co., Ltd.*
Xuchang China
rxsdliululu@126.com

*Abstract*—The diode rectifier unit-based high voltage direct current (DRU-HVDC) transmission with grid-forming (GFM) wind turbine is becoming a promising scheme for offshore wind farm(OWF) integration due to its high reliability and low cost. In this scheme, the AC network of the OWF and the DRU has completely different synchronization mechanisms and power flow characteristics from the traditional power system. To optimize the power flow and reduce the net loss, this paper carries out the power flow modeling and optimization analysis for the DRU-HVDC transmission system with grid-forming OWFs. The influence of the DRU and the GFM wind turbines on the power flow of the system is analyzed. On this basis, improved constraint conditions are proposed and an optimal power flow (OPF) method is established. This method can minimize the power loss by adjusting the reactive power output of each wind turbine and internal network frequency. Finally, based on MATLAB, this paper uses YALMIP toolkit and CPLEX mathematical solver to realize the programming solution of the OPF model proposed in this paper. The results show that the proposed optimization strategy can effectively reduce the power loss of the entire OWF and the transmission system with an optimization ratio of network losses exceeding 25.3%.

## I. Introduction

The modular multilevel converter-based high voltage direct current (MMC-HVDC) transmission system [1-2] is widely used to integrate remote offshore wind power into the onshore grid [3-4]. However, the large number of submodules in the offshore MMC station occupy a huge offshore platform area and thus greatly increase the construction cost. Consequently, there is an urgent need to develop a low-cost HVDC integration method for remote OWFs [5-6].

Compared with MMC-HVDC, the diode rectifier unit-based HVDC (DRU-HVDC) transmission system has a much smaller size, lighter weight, and lower loss. Therefore, the DRU-HVDC transmission is becoming a promising low-cost scheme to integrate offshore wind power. Due to the uncontrollable and unidirectional characteristics, the DRU cannot independently establish the offshore AC system or provide black-start power to the OWF. Conventional grid-following wind turbines (WTs) are not able to operate with DRU-HVDC transmission. To solve these problems, the grid-forming (GFM) control is proposed for wind turbines to establish the offshore AC system. In addition, energy storage devices are installed in some wind turbines to achieve black start [7-9].

Since the DRU-HVDC system is similar to a resistive load, the power flow characteristic between GFM WTs and the DRU-HVDC system is opposite compared to the traditional power system. The active power flow is dominated by the AC voltage amplitude, and the reactive power flow is related to the AC frequency. Currently, research on the optimization of the reactive power flow in such an AC system is scarce. Existing optimization analyses are mostly based on the MMC-HVDC system and grid-following (GFL) wind farms [10], where the system frequency characteristics are clear, the reactive power flow form is simple, and the modeling, as well as the solving, is easy. However, in the DRU-HVDC system with GFM WTs, the reactive power flow, reactive power constraints, and reactive power distribution are much more complex. Moreover, the diode reactive power characteristics have not been accurately analyzed yet, which makes it difficult to establish the system power flow model or to optimize it. To solve those problems, a comprehensive system power flow modeling is carried out in this paper, including the DRU-HVDC, GFM WTs, transformers, AC and DC cables, etc. Based on this model, the optimal power flow (OPF) analysis is completed, and an optimization strategy is proposed for WTs. The optimization objectives, including minimizing the system loss/maximizing the delivered active

This work is supported by the Project of IGCT-Based Novel Converter Development by China Electric Equipment Group Co., Ltd. (CEE-2022-B-01-01-006-XJ).

power, are realized by adjusting the reactive power output of each wind turbine and internal network frequency.

Finally, MATLAB is used as the solution platform, and the OPF problem proposed in this paper is solved based on the YALMIP toolkit and CPLEX mathematical planning problem solver. The proposed reactive power optimization strategy is verified by the optimal power flow results.

## II. System Topology and Control Structure Overview

### A. System Topology Overview

The reactive power self-synchronization control diagram based on GFM wind turbine is shown in Fig. 1, which includes a machine-side converter (MSC), a DC/DC energy storage converter, a grid-side converter (GSC), a DC bus and a diode uncontrolled rectifier unit.

### B. Control Structure Overview

The control of the GFM WT is proposed in [12]. As shown in Fig. 1, the control of the machine-side converter (MSC) remains the same with the conventional GFL control strategy. The modifications lie in the GSC: one is the introduction of an AC voltage middle-loop in the $d$-axis control loop to realize grid-forming characteristics; the other is the introduction of a distributed frequency-reactive power self-balancing control strategy in the $q$-axis control loop, which enables the OWF integration to have the capability of autonomous frequency construction and multi-turbine synchronization.

Through GFM wind turbines, the active-frequency self-synchronization mechanism brought by the swing equation of the synchronous machine is replaced by the dependence on the frequency-reactive power coupling control and the phase-locked loop (PLL) frequency tracking process.

In addition to its synchronization function, the PLL also determines the output frequency of the GSC by rotating the coordinate system, consequently defining the frequency of the offshore AC network within the wind farm. The mathematical expression for the PLL is shown in (1),

$$\omega = k_{Lp} u_{cfq} + k_{Li} \int u_{cfq} dt + \omega_0 \qquad (1)$$

where $k_{Lp}$ and $k_{Li}$ represent the proportional and integral coefficients of the PLL PI controller, respectively; $u_{cfq}$ is the $q$-axis voltage of the filtering capacitor $C_F$; $\omega_0$ is the rated frequency. Neglecting the damping resistor, the dynamic equation for $u_{cfq}$ is given by (2).

$$\begin{cases} C_F \dfrac{du_{cfq}}{dt} = i_{cfq} - \omega C_F u_{cfd} \\ i_{cfq} = i_{gcq} - i_{drq} \end{cases} \qquad (2)$$

In (2), $u_{cfd}$ and $u_{cfq}$ are the $d$-axis voltage and $q$-axis voltage of PCC; $i_{cfq}$, $i_{gcq}$, and $i_{drq}$ are filtering capacitor current, GSC output current, and wind farm output current, respectively. Considering that the output reactive power of GSC can be expressed as (3),

$$i_{gcq} = -Q_{wt}/1.5 u_{cfd} \qquad (3)$$

we can then acquire the coupling relationship between PLL frequency $\omega$ and GSC output reactive power $Q_{wt}$, which shows the basic principle of the frequency-reactive power coupling control mentioned above.

## III. Modeling of Optimal Power Flow Problem based on Frequency-Reactive Power Characteristics of OWF Integration

### A. Basic Optimal Flow Model

The OPF problem is the most common and fundamental optimization problem in power systems. On the premise of satisfying the physical constraints of power network such as Kirchhoff's law, line capacity constraints, and operational security constraints, the OPF problem aims to find an optimal power flow steady state operating point at which the objective functions of total system power generation cost, total network loss, etc. are optimized.

Referring to [11], several mathematical notations employed in the basic OPF model are introduced. Defining $N$ as the set of all nodes in the network, $E^+ = \{(i,j)|i<j$ and $i,j \in N\}$ as the set of directed branches, which represents a directed branch with two end nodes $i, j$ and positive direction $i$ to $j$, then the OPF model can be expressed in the form of the following (4) through (9).

$$\min_{S_i^g, V_i(\forall i \in N)} \sum_{\substack{i \in N \\ k=1,2}} c_{ki} \left( \Re \left( S_i^g \right) \right)^2 + c_{0i} \qquad (4)$$

$$s.t. \ \underline{V} \leq |V_i| \leq \overline{V}, \ \forall i \in N \qquad (5)$$

$$\underline{S}_i^g \leq S_i^g \leq \overline{S}_i^g, \ \forall i \in N \qquad (6)$$

$$|S_{ij}| \leq \overline{S}_{ij}, \ \forall (i,j) \in E^+ \cup E^- \qquad (7)$$

$$S_i^g - S_i^d = \sum_{(i,j) \in E^+ \cup E^-} S_{ij}, \ \forall i \in N \qquad (8)$$

$$S_{ij} = Y_{ij}^* V_i V_i^* - Y_{ij}^* V_i V_j^*, \ \forall (i,j) \in E^+ \cup E^- \qquad (9)$$

In the above basic OPF model, $V_i$ and $S_i^g$ represent the voltage and injected power of node $i$ respectively (both are complex variables); $S_i^d$ represents the load of node $I$; $S_{ij}$ represents the complex power flowing from node $i$ into line $(i, j)$ (note that due to the line losses, $S_{ij}+S_{ji}\neq 0$); $c_{xi}$ ($x$=0,1,2) represents the coefficients of generation cost in the objective function. The objective function is shown in (4), where $\Re(\cdot)$ is a function that takes the real part of the complex variable, this is because the classical OPF objective function (e.g., generation cost) is generally only related to the injected active power of nodes. The constraints on the node voltage magnitude, node injected power, and line transmission power are shown in (5) to (7) respectively. KCL and KVL laws are shown in (8) and (9) respectively.

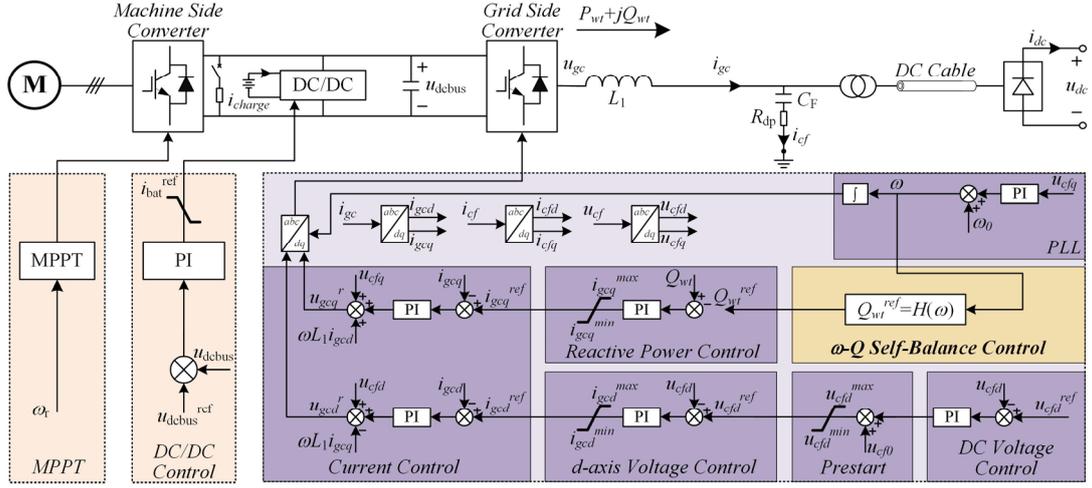

Fig. 1. The reactive power synchronization control diagram based on GFM wind turbines

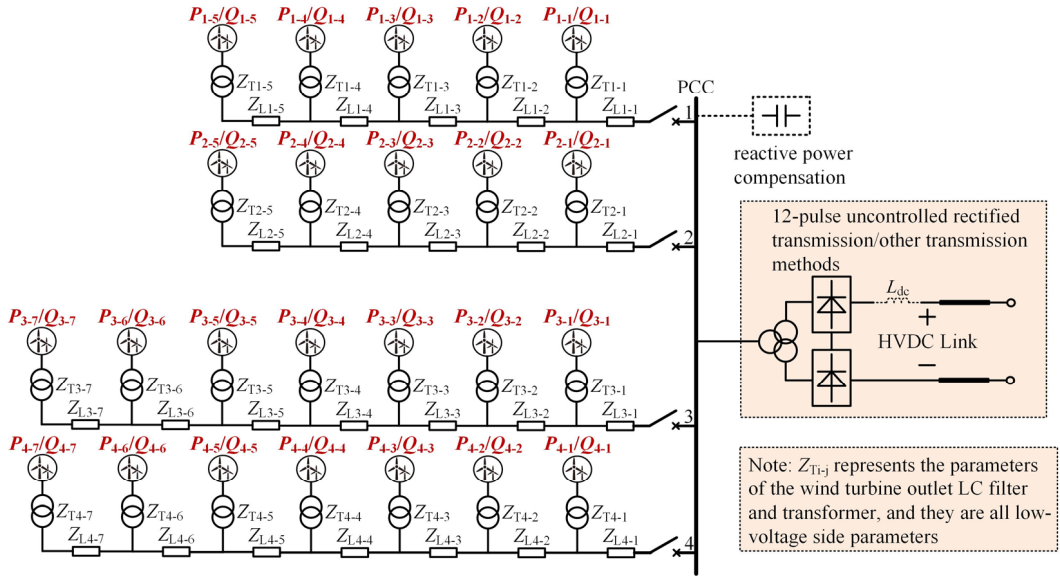

Fig. 2. The reactive power synchronization control diagram based on GFM wind turbines

## B. Optimal Power Flow Model Considering Frequency-Reactive Power Characteristics of OWF Integration

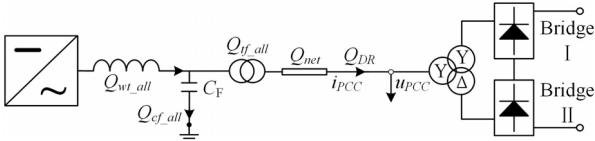

Fig. 3. Reactive power distribution diagram of wind farm

In different practical application scenarios, the constraints of the basic OPF model described in Section III.A may be increased. As for the model studied in this paper, the reactive power flow distribution and frequency characteristics of the system need to be considered more comprehensively. The static operating point of the wind farm integration frequency depends on the frequency-reactive power output characteristics of the offshore wind turbine and the total reactive power demand of the integration system. The distribution of reactive power supply and demand in the wind farm integration is shown in Fig. 3, where $Q_{wt\_all}$ represents the total reactive power output of the GSC; $Q_{cf\_all}$ represents the total reactive power consumption of the filter capacitor $C_F$; $Q_{tf\_all}$ represents the total reactive power consumption of the turbine transformer; $Q_{net}$ represents the reactive power consumed by the collection network; $Q_{DR}$ represents the reactive power consumed by the diode rectifiers. Besides, $u_{pcc}$ and $i_{pcc}$ represent the phase voltage peak and the phase current peak at the point of common coupling (PCC) respectively.

Assuming that all wind turbines have the same parameters and work in the same operating condition, the wind farm reactive power balance equation can be expressed as (10),

$$\begin{cases} Q_{wt\_all}(\omega) = Q_{farm}(P_{farm}, \omega) \\ Q_{farm} = Q_{cf\_all} + Q_{tf\_all} + Q_{net} + Q_{DR} \end{cases} \quad (10)$$

where $P_{farm}$ represents the active power output of the wind farm and $Q_{farm}$ represents the reactive power consumed by the whole system. $Q_{DR}$, $Q_{net}$, $Q_{tf\_all}$ and $Q_{cf\_all}$ can be calculated by the following (11) to (15) respectively.

$$Q_{DR} = P_{farm} \tan(\varphi) \quad (11)$$

$$\varphi = \arctan\left(\frac{2\mu - \sin(2\mu)}{1 - \cos(2\mu)}\right) \quad (12)$$

$$\begin{cases} Q_{net} = 1.5 i_{pcc}^2 \omega L_{net} - 1.5 u_{pcc}^2 \omega C_{net} \\ i_{pcc} = P_{farm}/1.5 u_{pcc} \cos(\varphi) \end{cases} \quad (13)$$

$$Q_{tf\_all} = 1.5 N_{\omega t}\left(i_{pcc} n_{tf}/N_{\omega t}\right) \omega L_{tf} \quad (14)$$

$$Q_{cf\_all} = -1.5\left(u_{pcc}/n_{tf}\right)^2 N_{\omega t} \omega C_F \quad (15)$$

In (11), $\varphi$ represents the power factor angle at PCC, with its calculation formula referenced in (12); $N_{wt}$ represents the number of wind turbines; $L_{net}$ and $C_{net}$ represent the centralized equivalent values of inductance and capacitance in the collection network respectively; $n_{tf}$ and $L_{tf}$ represent the turbine transformer ratio and low-voltage side equivalent leakage inductance. For the system described in this paper, based on (11) to (15), a power flow model is established and the normalized frequency range of the OWF integration is set as (0.9, 1.1). When the active power output of wind farm uniformly increases from 0.01p.u. to 0.8p.u., the curve of reactive power consumption in the OWF integration is shown in Fig. 4. The horizontal axis in the figure represents the normalized value of $\omega$, and the vertical axis represents the total reactive power output demand. From this figure, it is obvious that for a fixed wind farm active power output, the relationship between the total reactive power demand and the normalized frequency is approximately linear. What's more, it exhibits a negative correlation at lower active power output and a positive correlation at higher active power output.

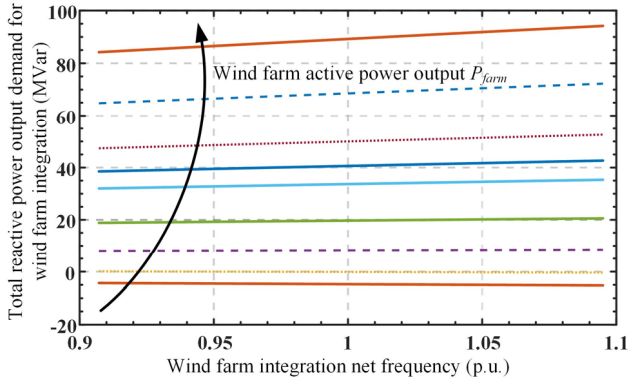

Fig. 4. Static operation point of wind farm frequency (the arrow indicates the increasing active power)

As shown in Fig. 4, when the internal network frequency varies within the desired range of $(\omega_{min\_h}, \omega_{max\_h})$, by selecting control parameters correctly, the static operating point of the system can be made to reach any position on any curves in the graph. The calculation approach for control parameters is described by (16) to (20) below,

$$k_H = \frac{(Q_{farm\_max\_h} - Q_{farm\_min\_h})}{N_{\omega t} \Delta \omega_{max}} \quad (16)$$

$$Q_0 = \frac{Q_{farm\_max\_h}}{N_{\omega t}} - k_H\left(\omega_{max\_h} - \omega_0\right) \quad (17)$$

$$\Delta \omega_{max} = \omega_{max\_h} - \omega_{min\_h} \quad (18)$$

$$\omega_{min\_h} \leq \omega_0 \leq \omega_{max\_h} \quad (19)$$

$$Q_{farm\_max\_h} = Q_{farm}\left(\omega_{max\_h}, P_{farm\_max}\right) \quad (20)$$

where $Q_{farm\_min\_h}$ represents the reactive power consumed by the wind farm when the internal network frequency is $\omega_{min\_h}$ and the wind farm's active power output is 0. $Q_{farm\_max\_h}$ represents the reactive power consumed by the wind farm when the internal network frequency is $\omega_{max\_h}$ and the wind farm is operating at full capacity. $P_{farm\_max}$ represents the wind farm's maximum active power output.

### C. Model Convexification and Relaxation

Combining Sections III.A and III.B, it can be observed that the non-convexity of the OPF model constructed in (4) to (15) arises from two aspects: the first aspect is the strong non-linearity introduced by (11) and (12); the second aspect originates from the quadratic terms in (9) and the absolute value lower bound in (5), i.e., the product of the voltage complex variables and the logic relationship that the absolute value should exceed the lower bound. To solve the non-convex problem in the model and obtain the optimal solution within polynomial time, further adjustments to the model are required.

For the first aspect, although these two equations exhibit strong non-linearity, as can be seen in Fig. 4, the relationship between the desired reactive power output and the frequency is actually quite close to a linear relationship. Meanwhile, under the given parameter conditions, numerical analysis of the non-linear relationship introduced by (11) and (12), namely the relationship between $Q_{DR}$ and the transformer's low-voltage side voltage $u_P$, can be obtained, as shown in Fig. 5.

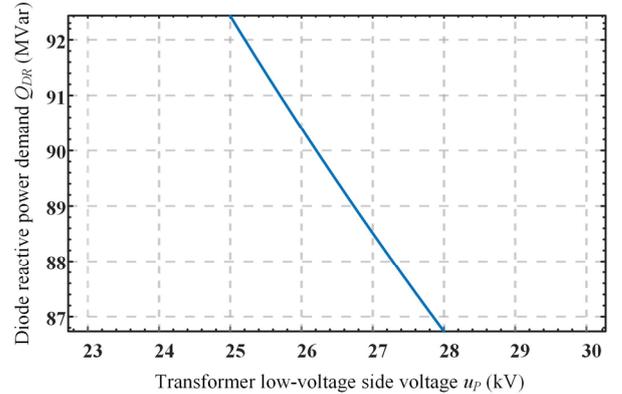

Fig. 5. Nonlinear reactive power flow characteristics introduced by diode rectifier

From this figure, it can be seen that although (11) and (12) introduce non-linear characteristics, the numerical solution results of the non-linear relationships within the dynamic operating range corresponding to the example parameter conditions exhibit strong linearity. Therefore, these non-linear characteristics brought by (11) and (12) can be approximated using the Taylor expansion method, and only consider the first two terms of the Taylor expansion for model development. This leads to the reactive power demand constraint as shown in (21) below,

$$Q_{\omega t\_all}(\omega, P_1, P_2, \ldots, P_k) = k_1\left(\sum P_i\right)^2 \omega + k_2\omega + k_3\sum P_i \quad (21)$$

where $P_i$ ($i = 1,2,3,\ldots,k$) represents active power for each wind turbine; $k_j$ ($j = 1,2,3$) represents coefficients. However, in the actual operation of a wind farm, to maximize the utilization of wind energy, the power control of converter often employs the maximum power point tracking (MPPT) algorithm. As a result, the active power generated by the wind turbines changes along with the time. We simplify the problem by assuming that the active power of wind turbines will remain constant value determined by environmental factors at different time intervals. This simplifies (21) to the linear reactive power constraint shown in (22) for each time interval. Without loss of generality, assuming the active power of all wind turbines are the same and we focus on the active power of a single wind turbine. Taking the case of a single wind turbine with an active power output of 4.54 MW as an example (a typical value within the 24-hour dynamic active power output change), the calculated coefficients for the linear reactive power demand constraint are as shown in (23).

$$Q_{\omega t\_all}(\omega) = d_1\omega + d_2 \quad (22)$$

$$\begin{cases} d_1 = 0.0127 \\ d_2 = 36.6340 \end{cases} \quad (23)$$

As for the non-convex factors in the second aspect, the common approach is to introduce a new variable $W$ to represent the product term of the voltage complex variables.

$$V_i V_j^* = W_{ij}, \quad i,j \in N \quad (24)$$

Then rewriting (5) and (9) based on (24) yields:

$$\underline{V}^2 \leq W_{ii} \leq \overline{V}^2, \quad \forall i \in N \quad (25)$$

$$S_{ij} = Y_{ij}^* W_{ii} - Y_{ij}^* W_{ij}, \quad \forall (i,j) \in E^+ \cup E^- \quad (26)$$

Denote the voltage column vector to be solved as $V$, then $W$ is a rank-1 positive semidefinite matrix variable. Once the value of $W$ is obtained, the unique solution vector can be recovered from it. The SDP relaxation (also known as Shor relaxation) involves the abandonment of the challenging non-convex constraint of rank($W$)=1 in the model, thereby relaxing the OPF problem into an SDP programming problem. To further simplify this problem, considering the semi-definite equivalence conditions for $W$ regarding principal minors, ignoring the inequalities involving principal minor of the third order and above, yields a specific class of SDP relaxation known as second-order cone programming (SOCP) relaxation. The standard second-order cone constraints are as follows, and the SOCP relaxation is the convex relaxation strategy employed in this paper.

$$\left\| \begin{matrix} 2W_{ij} \\ W_{ii} - W_{jj} \end{matrix} \right\| \leq W_{ii} + W_{jj} \quad (27)$$

## IV. MODEL SOLVING AND CONCLUSIONS

### A. Overview of Model Solving

Based on the system topology shown in Fig. 2, the model is solved using feeder parameters, LCL filter parameters, wind turbine parameters, and transformer parameters, which all comes from the real project.

The model is solved using the YALMIP toolkit and the CPLEX mathematical programming solver. Based on the dynamic active power output over 24 hours from a specific wind farm, the frequency-reactive power optimization strategy for OWF integration based on GFM wind turbines is implemented and verified.

### B. Model Solving and Conclusions

The solving results of the OPF model established in this paper are presented in Fig. 6, which demonstrate the internal network frequency value of optimization result as well as a comparison of network loss between the original strategy and proposed optimal strategy. From Fig. 6 (a), which shows the values of dynamic frequency optimization in 24 hours, it can be observed that the frequency variation range is (0.9995,1.0001), which is relatively small. Meanwhile, as is shown in Fig. 6 (b), the results of the 24-hour dynamic network loss optimization reveal that the proposed reactive power output optimization control strategy is effective in reducing system network losses when the wind turbine's active power output is relatively high. The optimization ratio of network loss remains stable at 25.3% or higher from 1h to 14h of the day. However, when the wind turbine's active power output is relatively low, due to the small network losses, the modeling approximations and simplifications mentioned above, the network losses either approach or slightly exceed those incurred under the original reactive power control strategy.

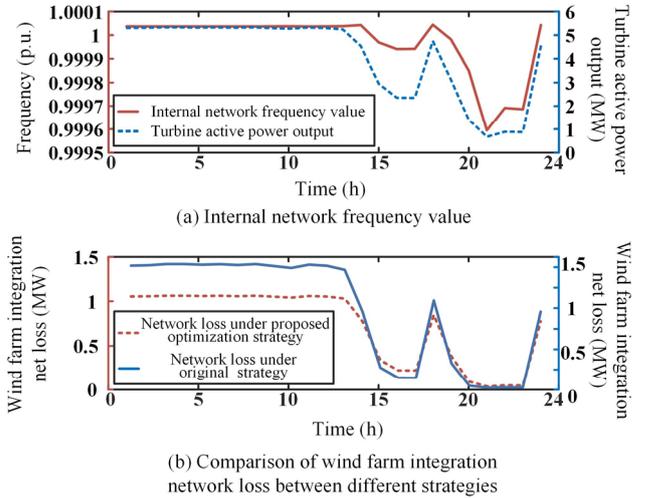

(a) Internal network frequency value

(b) Comparison of wind farm integration network loss between different strategies

Fig. 6. Optimal frequency value of wind farm integration and comparison of network loss optimization results

The comparison of the reactive power output by the wind turbines at different time intervals is shown in Fig. 7. In Fig. 7 (a), the horizontal axes represent the time division of a day in 24 hours and the unit numbers of 12 wind turbines from two feeder clusters respectively. The vertical axis represents the reactive power output of the wind turbines.

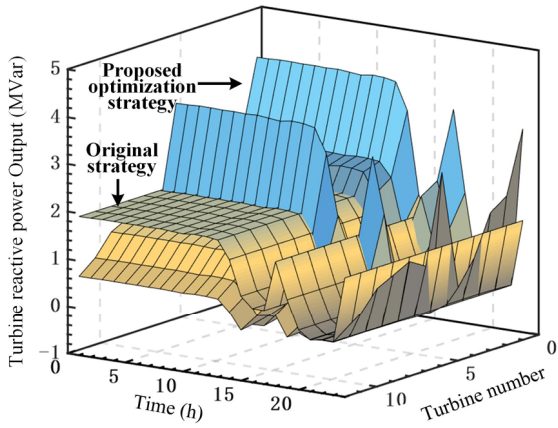

(a) Comparison of reactive power output of wind farms at different time Intervals – 3D view

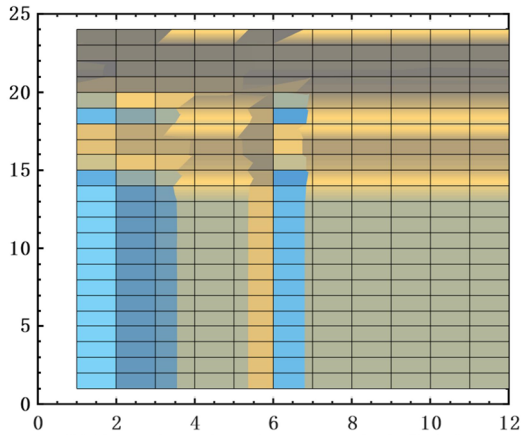

(b) Comparison of reactive power output of wind farms at different time intervals - 2D view

Fig. 7. Comparison of reactive power output of wind farm in different time periods

The figure displays the reactive power output distribution of the wind turbines under two different control strategies at various time intervals, where the fluctuating curve with prominent blue peaks represents the three-dimensional waveform of reactive power output under the proposed frequency-reactive power optimization control strategy. In contrast, the smoother curve represents the waveform of reactive power output under the original uniform reactive power output strategy. Furthermore, for a clearer comparison of the magnitude of the reactive power output, a two-dimensional view is used to analyze and contrast the distribution of reactive power output.

## V. CONCLUSIONS

Using GFM WTs and DRU-HVDC transmission can realize low-cost grid integration of remote OWFs. This paper focus on the power interaction between GFM WTs and DRU in offshore AC system. A detailed model is established for power flow analysis and optimization. Based on the basic OPF model, improved optimization constraints, including reactive power demand and frequency stability, are considered. A frequency-reactive power optimization control strategy is conducted by adjusting the reactive power output of each wind turbine and the internal network frequency with the objective of minimizing system network losses/maximizing the deliverer active power. The simulation results shows that the proposed optimization strategy can effectively reduce network losses for the offshore AC system. This effect is particularly significant when the active power output of WTs is relatively high (50% to 70% load capacity), with an optimization ratio of network losses exceeding 25.3%.